\title{Values of coefficients of cyclotomic polynomials II}
\author{Chun-Gang Ji, Wei-Ping Li and Pieter Moree}
\documentclass[12pt]{article}
\usepackage{amssymb, latexsym, amsfonts}
\def\@ptsize{2}
\newtheorem{Thm}{Theorem}

\newtheorem{Lem}{Lemma}

\newtheorem{cor}{Corollary}

\newcommand{\qed}{\hfill $\Box$}

\begin{document}
\date{}
\maketitle
{\def\thefootnote{}
\footnote{
The first author is partially supported by the Grant No.10771103 from NNSF of China.}

\begin{abstract} \noindent
Let $a(n, k)$ be the $k$th coefficient of the $n$th cyclotomic
 polynomial. In part I it was proved that  $\{ a(mn, k)\  | \ n\ge 1,\ k\ge 0 \}={\mathbb Z}$, in case $m$ is a prime power. In this paper
 we show that the result also holds true in case $m$ is an arbitrary positive
 integer.
\end{abstract}
\section{Introduction}
Let $\Phi_n(x)=\sum_{k=0}^{\varphi(n)}a(n,k)x^k$ be the $n$th cyclotomic polynomial.
The rational function $1/\Phi_n(x)$ has a Taylor series around $x=0$ given by
$${1\over \Phi_n(x)}=\sum_{k=0}^{\infty}c(n,k)x^k,$$
where it can be shown that the $c(n,k)$ are also integers. 
It turns out that usually the coefficients $a(n,k)$ and $c(n,k)$ are quite small
in absolute value, for example for $n<105$ it is well-known that $|a(n,k)|\le 1$ and for $n<561$ we have
$|c(n,k)|\le 1$ (by \cite[Lemma 12]{Pieter}).\\
\indent The purpose of this note
is to show that although so often the coefficients $a(n,k)$ and $c(n,k)$ are small, they
assume every integer value, even when we require $n$ to be a multiple of an arbitrary natural
number $m$.
\begin{Thm} Let $m\ge 1$ be an integer. Put
$S(m)=\{ a(mn, k) | n\ge 1,~k\ge 0 \}$ and
$R(m)=\{ c(mn, k) | n\ge 1,~k\ge 0\}$. Then
$S(m)=\mathbb Z$ and $R(m)=\mathbb Z$.
\end{Thm}
Schur poved in 1931 (in a letter to E. Landau) that $S(1)$ is not a finite set. In 1987 Suzuki \cite{Su} proved
that $S(1)=\mathbb Z$. Recently the first two authors \cite{jili} proved
that $S(p^e)=\mathbb Z$ with $p^e$ a prime power.\\
\indent The fact that every integer already occurs as a coefficient of 
$\Phi_{pqr}(x)$ with $p$, $q$ and $r$ odd
primes is implicit in Bachman \cite{Ba1}. The third author established this
result for the reciprocal cyclotomic polynomials $1/\Phi_{pqr}(x)$, see Moree \cite{Pieter}.
This result implies that $R(1)=\mathbb Z$.

\section{Some lemmas}
Since 
\begin{equation}
\label{flauw}
x^n-1=\prod\limits_{d|n}\Phi_{d}(x), 
\end{equation}
we have by the M\"{o}bius
inversion formula,
$\Phi_n(x)=\prod_{d|n}(x^d-1)^{\mu(\frac{n}{d})}$,
where $\mu$ denotes the M\"{o}bius function.

On using that $\sum_{d|n}\mu(d)=0$ if $n>1$, it is seen that, for
$n>1$,
$$\Phi_n(x)=\prod_{d|n}(x^d-1)^{\mu(\frac{n}{d})}=
(-1)^{\sum_{d|n}\mu(\frac{n}{d})}\prod_{d|n}(1-x^d)^{\mu(\frac{n}{d})}
=\prod_{d|n}(1-x^d)^{\mu(\frac{n}{d})}.
$$
(Thus for $n>1$, the polynomial $\Phi_n(x)$ is self-reciprocal.)
\begin{Lem}
\label{conge}
The coefficient $c(n,k)$ is an integer whose values only depends on the
congruence class of $k$ modulo $n$.
\end{Lem}
{\it Proof}. Let us first consider
$$\Psi_n(x):={x^n-1\over \Phi_n(x)}.$$
By (\ref{flauw}) we have that $\Psi_n(x)=\prod_{d<n,~d|n}\Phi_d(x)$ and
thus its coefficients are integers. The degree of $\Psi_n(x)$ is 
$n-\varphi(n)$, with $\varphi$ Euler's totient function. We infer that, for $|x|<1$,
$${1\over \Phi_n(x)}=-\Psi_n(x)(1+x^n+x^{2n}+\cdots )$$
Since $n>n-\varphi(n)$, the proof is completed. \qed\\

\noindent Let $\kappa(m)=\prod_{p|m}p$ denote the squarefree kernel of $m$, that is
the largest squarefree divisor of $m$.
\begin{Lem}
\label{blue} Let $p$ be a prime.
For $l,m\geq 1$ we have $S(p^lm)=S(pm)$ and $R(p^lm)=R(pm)$.
\end{Lem}
\begin{cor}
We have $S(m)=S(\kappa(m))$ and $R(m)=R(\kappa(m))$.
\end{cor}
{\it Proof of Lemma} \ref{blue}. It is easy to prove,
see e.g. Thangadurai \cite{Th}, that if $p$ is prime and $p|n $, then
\begin{equation}
\label{blo}
\Phi_{pn}(x)=\Phi_n(x^p).
\end{equation}
Using this we deduce 
that
$\Phi_{p^2m}(x)=\Phi_{pm}(x^p)$ and thus 
$a(pm,1)=0$ and hence $0\in S(pm)$. On repeatedly applying
(\ref{blo}) we can easily infer that
$\Phi_{p^lmn}(x)=\Phi_{pmn}(x^{p^{l-1}})$ for any $l\geq 1$, so
$$a(p^lmn,k)=\cases{a(pmn,{k\over p^{l-1}}) & if $p^{l-1}|k$;\cr
0 & otherwise.}$$ 
This together with
$0\in S(pm)$ and the trivial inclusion
$S(p^lm)\subseteq S(pm)$ shows that $S(p^lm)=S(pm)$.\\
\indent The proof that $R(p^lm)=R(pm)$ is completely analogous.
Here we use that if $p|n$, then $\Psi_{pn}(x)=\Psi_n(x^p)$,
which is immediate from (\ref{blo}) and the definition of $\Psi_n(x)$. \qed

\begin{Lem} {\rm (Quantitative form of Dirichlet's theorem.)} Let $a$
and $m$ be coprime natural numbers and let $\pi(x; m, a)$ denote the
number of primes $p\le x$ that satisfy $p\equiv a({\rm mod~}m)$. Then,
as $x$ tends to infinity,
$$\pi(x; m, a)\sim \frac{x}{\varphi(m)\log x}.$$
\end{Lem}
\begin{cor} Given $m,t\geq 1$ and any real number $r>1$, there exists a
constant $N_0(t, m,r)$ such that for every $n>N_0(t, m,r)$ the interval
$(n, rn)$ contains at least $t$ primes $p\equiv 1({\rm mod~}m)$.
\end{cor}

\section{The proof of Theorem 1}
We first prove that $S(m)=\mathbb Z$.
Since $S(m)=S(\kappa(m))$, we may assume that $m$ is squarefree.
We may also assume that $m>1$.
Suppose that $n>N_0(t,m,{15\over 8})$. Then there exist
primes $p_1$, $p_2$, $\cdots$, $p_t$ such that
$$n<p_1<p_2<\cdots <p_t<{15\over 8}n{\rm ~and~}
p_j\equiv 1 ({\rm mod~} m), \quad j=1,2, \cdots, t.$$
Hence $p_t<2p_1$.

Let $q$ be any prime exceeding $2p_1$ and put
$$m_1=\cases{p_1p_2\cdots p_tq & if $t$ is even;\cr
p_1p_2\cdots p_t & otherwise.}$$ 
Note that $m$ and $m_1$ are coprime and that $\mu(m_1)=-1$. Using
these observations we conclude that
\begin{eqnarray}
\label{bluebird}
\Phi_{m_1m}(x)&\equiv &\prod_{d|m_1m,~d<2p_1}(1-x^d)^{\mu({m_1m\over
d})} ~({\rm mod~}x^{2p_1})\cr
&\equiv &\prod_{d|m}(1-x^d)^{\mu({m\over
d})\mu(m_1)}\prod_{j=1}^t(1-x^{p_j})^{\mu({m_1m\over p_j})} ~({\rm mod~}x^{2p_1}) \cr
&\equiv & \Phi_{m}(x)^{\mu(m_1)}
\prod_{j=1}^t(1-x^{p_j})^{-\mu(m_1m)} ~({\rm mod~}x^{2p_1}).\cr
&\equiv & {1\over \Phi_{m}(x)}
\prod_{j=1}^t(1-x^{p_j})^{\mu(m)} ~({\rm mod~}x^{2p_1}).\cr
&\equiv & {1\over \Phi_{m}(x)}\Big(1-\mu(m)(x^{p_1}+\ldots+x^{p_t})\Big)~({\rm mod~}x^{2p_1}).
\end{eqnarray} 
From (\ref{bluebird}) it follows that, if $p_t\le k<2p_1$,
$$a(m_1m,k)=c(m,k)-\mu(m)\sum_{j=1}^tc(m,k-p_j).$$ By
Lemma \ref{conge} we have $c(m,k-p_j)=c(m,k-1)$. Thus we find that
\begin{equation}
\label{boehoe}
a(m_1m,k)=c(m,k)-\mu(m)tc(m,k-1) {\rm ~with~}p_t\le k<2p_1.
\end{equation}
\indent We
consider two cases ($\mu(m)=1$, respectively $\mu(m)=-1$).\\
{\bf Case 1}. $\mu(m)=1$. In this case $m$ has at least two prime divisors.
Let $q_1<q_2$ be the smallest two prime divisors of $m$. Here
we also require that $n\ge 8q_2$. This ensures that $p_t+q_2<2p_1$.
Note that
\begin{eqnarray}
\label{eagle}
{1\over \Phi_m(x)} & \equiv & {(1-x^{q_1})(1-x^{q_2})\over 1-x} ~({\rm mod~}x^{q_2+2})\cr
& \equiv & 1+x+x^2+\ldots+x^{q_1-1}-x^{q_2}-x^{q_2+1} ~({\rm mod~}x^{q_2+2}).
\end{eqnarray}
Thus $c(m,k)=1$ if $k\equiv \beta ({\rm mod~}m)$ with $\beta\in \{1,2\}$
and $c(m,k)=-1$ if $k\equiv \beta ({\rm mod~}m)$ with $\beta\in \{q_2,q_2+1\}$.
This in combination with (\ref{boehoe}) shows that $a(m_1m,p_t+1)=1-t$
and $a(m_1m,p_t+q_2)=t-1$. Since $\{1-t,t-1\ | \ t\ge 1\}=\mathbb Z$ the result
follows in this case.\\
{\bf Case 2}. $\mu(m)=-1$. Here we notice that
$${1\over \Phi_m(x)}\equiv \cases{1-x ~({\rm mod~}x^3) & if $2\nmid m$;\cr
1-x+x^2 ~({\rm mod~}x^3) & otherwise.}$$
Using this we find that $a(m_1m,p_t)=-1+t$. 
Furthermore, $a(m_1m,p_t+1)=-t$ in case $m$ is
odd and $a(m_1m,p_t+1)=1-t$ otherwise. Since
$\{-1+t,-t\ | \ t\ge 1\}=\mathbb Z$ and $\{-1+t,1-t\ | \ t\ge 1\}=\mathbb Z$,
it follows that also $S(m)=\mathbb Z$ in this case.

It remains to show that $R(m)=\mathbb Z$. As before we may assume that $m$ is
squarefree (by Corollary 1) and that $m>1$ (by Theorem 8 of Moree \cite{Pieter}).\\
\indent Let $q$ be any prime exceeding $2p_1$ and put
$$\overline{m}_1=\cases{p_1p_2\cdots p_t & if $t$ is even;\cr
p_1p_2\cdots p_tq & otherwise.}$$
Note that $\mu(\overline{m}_1)=1$. Reasoning as in the derivation of (\ref{bluebird})
we obtain
$${1\over \Phi_{\bar{m}_1m}(x)}\equiv {1\over \Phi_{m}(x)}\Big(1-\mu(m)(x^{p_1}+\ldots+x^{p_t})\Big)~
({\rm mod~}x^{2p_1})$$
and from this
$c(\bar{m}_1m,k)=a(m_1m,k)$ for $k<2p_1$. Reasoning as in the proof of $S(m)=\mathbb Z$, the
proof is then completed. \qed\\

\noindent {\tt Remark 1}. If one specializes the above proof to the case $m=p^e$, a proof
a little easier than that given in part I \cite{jili} is obtained, since it does not involve
a case distinction between $m$ is odd and $m$ is even as made in part I. This
is a consequence of working modulo $x^{2p_1}$, rather than modulo $x^{2p_1+1}$.\\

\noindent {\tt Remark 2}. The fraction $15/8$ in the proof can be replaced by $2-\epsilon$,
with $0<\epsilon<1$ arbitrary. One then requires that $n>N_0(t,m,2-\epsilon)$ and
in case $\mu(m)=1$ in addition that $n\ge q_2/\epsilon$. 

\vfil\eject

\medskip\noindent {\footnotesize Department of Mathematics,\\
Nanjing Normal University\\
Nanjing 210097, P.R. China\\
e-mail: {\tt cgji@njnu.edu.cn}}

\medskip\noindent {\footnotesize Rugao Normal College,\\
Rugao 226500, Jiangsu, P.R. China\\
e-mail: {\tt lwpeace@sina.com}}

\medskip\noindent {\footnotesize Max-Planck-Institut f\"ur Mathematik,\\
Vivatsgasse 7, D-53111 Bonn, Germany.\\
e-mail: {\tt moree@mpim-bonn.mpg.de}}

\end{document}